\documentclass[12pt]{amsart}
\usepackage{amscd,amssymb}
\usepackage[arrow,matrix]{xy}

\theoremstyle{plain}
\newtheorem{theorem}[subsection]{Theorem}
\newtheorem{lemma}[subsection]{Lemma}
\newtheorem{prop}[subsection]{Proposition}
\newtheorem{corollary}[subsection]{Corollary}

\theoremstyle{definition}
\newtheorem{remark}[subsection]{Remark}

\numberwithin{equation}{section}
\newcommand{\HC}{{\mathcal H}}
\newcommand{\F}{{\mathcal F}}
\newcommand{\G}{{\mathcal G}}

\newcommand{\LL}{{\mathcal L}}

\newcommand{\Z}{\mathbb{Z}}
\newcommand{\Q}{\mathbb{Q}}

\newcommand{\C}{\mathbb{C}}

\newcommand{\PP}{\mathbb{P}}
\newcommand{\HH}{\mathbb{H}}

\begin{document}
\date{}

\title [ Alexander Invariants and Transversality ]
{  Alexander Invariants  and Transversality    }

\author[Alexandru Dimca]{Alexandru Dimca}
\address{  Laboratoire J.A. Dieudonn\'e, UMR du CNRS 6621,
                 Universit\'e de Nice-Sophia-Antipolis,
                 Parc Valrose,
                 06108 Nice Cedex 02,
                 FRANCE.}
\email
{dimca@math.unice.fr}



\begin{abstract}
We show that some of the main results in Laurentiu Maxim's paper \cite{M} can be obtained (even in a slightly more general setting) using the theory of perverse sheaves of finite rank over $\Q$ as described for instance in
author's recent book  \cite{D}.

\end{abstract}

\maketitle

\section{The main results}

Let $X \subset \C^{n+1}$ with $n>1$ be a reduced hypersurface given by an equation $f=0$.
We say that $f: \C^{n+1} \to \C$ is $\infty$-transversal if the projective closure $V$ of $X$ in $ \PP^{n+1}$ is transversal in the stratified sense to the hyperplane at infinity $H=
\PP^{n+1}  \setminus  \C^{n+1}$.

  Consider the affine complement $M_X= \C^{n+1} \setminus X$, and denote by $M_X ^c$ its infinite cyclic covering corresponding to the morphism
$$f_{\sharp}:\pi _1(M_X) \to \pi _1(\C^*)=\Z.$$

Then, for any positive integer $m$, the homology group $H_m(M_X ^c,\C)$, regarded as a module over the principal ring $\Lambda =\C[t,t^{-1}]$ 
as in \cite{Lib}, \cite{DN},
is the $m$-th
 Alexander module of the hypersurface $X$. When this module is torsion, we denote by 
$\Delta _m(t)$ the corresponding  $m$-th  Alexander polynomial of $X$.

The aim of this note is to give an alternative proof for generalizations of Corollary 3.8 and Theorem 4.1, as well as for
 a reformulation of Theorem 4.2 in \cite{M}. With the above notation, the first two results in \cite{M} can be stated as follows.

\begin{theorem}  \label{unuth}
Assume that  $f: \C^{n+1} \to \C$ is $\infty$-transversal.
Then the  Alexander modules  $H_m(M_X ^c,\C)$   of the hypersurface $X$ are torsion for $m<n+1$.

\end{theorem}

Since $M_X ^c$ is an $(n+1)$-dimensional CW complex, one has $H_m(M_X ^c,\C)=0$  for $m>n+1$, while
  $H_{n+1}(M_X ^c,\C)$ is free. In this sense, the above result is optimal.

\begin{theorem}  \label{doith}
Assume that  $f: \C^{n+1} \to \C$ is $\infty$-transversal.
Let $ \lambda \in \C^{*}$ be such that $\lambda  ^d \ne 1$, where $d$ is the degree of $V$.
Then $ \lambda $ is not a root of the  Alexander polynomials $\Delta _m(t)$ for $m<n+1$.
\end{theorem}


Recall that the construction of the Alexander modules and polynomials was generalized in the obvious way
in  \cite{DN} to the case when $ \C^{n+1}$ is replaced by a smooth affine variety $U$.

Let $W'=W_0' \cup ... \cup W_m'$ be a hypersurface arrangement in $ \PP^{N}$ for $N>1$. Denote by
$M_{W'}= \PP^{N} \setminus W'$ the corresponding complement.
 Let $d_j$ denote the degree of $W_j'$ and let $g_j=0$ be a reduced defining equation for $W_j'$ in 
$ \PP^{N}$. Let $Z \subset \PP^{N}$ be a smooth complete intersection of dimension $n+1>1$ and set $W_j= W_j' \cap Z$ for $j=0,...,m$ be the corresponding hypersurfaces in $Z$ considered as subscheme
defined by the principal ideals generated by the $g_j$'s.

 Let $W=W_0 \cup ... \cup W_m$ denote
the corresponding hypersurface arrangement in $Z$.
We assume troughout in this paper that the following hold.

\bigskip

\noindent ({\bf H1}) All the hypersurfaces $W_j$ are distinct, reduced and irreducible; moreover  $W_0$ is smooth. 

\bigskip

\noindent ({\bf H2}   ) The hypersurface   $W_0$ is transverse in the stratified sense to $V=W_1 \cup ... \cup W_m$.

\bigskip

 The complement $U= Z \setminus W_0$ is a smooth affine variety and we consider the hypersurface $X= U \cap V$ in $U$ and its complement
 $M_W = Z \setminus W = U \setminus X $.

\begin{theorem}  \label{treith}

Assume that $d_0$ divides the sum $\sum_{j=1,m}d_j$, say $d d_0=\sum_{j=1,m}d_j$. Then 
one has the following.

\noindent (i) The function
$f:U \to \C$ given by
$$f(x)= \frac{g_1(x)...g_m(x)}{g_0(x)^d}$$
is a well-defined regular function on $U$ whose generic fiber $F$ is connected.  

\noindent (ii) The above Theorems \ref{unuth} and \ref{doith} hold for the Alexander modules and the Alexander polynomials associated to the hypersurface $X$ in $U$.
\end{theorem}

Note that we need the connectedness of $F$ since this is one of the general assumptions made in
\cite{DN}. The next result says roughly that an $\infty$-transversal polynomial behaves as a homogeneous polynomial up-to (co)homology of degree $n-1$. In these degrees, the determination of the  Alexander polynomial of $X$ in $U$ is reduced to the simpler problem of computing a monodromy operator.

\begin{corollary}  \label{top}
With the assumption in Theorem  \ref{treith}, the following hold.

\noindent (i) Let $\iota:\C^* \to\C$ be the inclusion. Then, for each  $m<n$ there is a local system $\LL_m$ on  $\C^*$ such that 
$$R^mf_*\C_U=\iota   _!\LL_m.$$
 In particular, for each $m<n$, the monodromy operators of $f$ at the origin $T_0^m$ and at infinity 
$T_{\infty}^m$ acting on
$H^m(F,\C)$  coincide and the above local system $\LL_m$ is precisely the local system corresponding to this automorphism of $H^m(F,\C)$.

\noindent (ii) For $m<n$, there is an isomorphism $H^m(F,\C) \to H^m(M_W^c,\C)$ which is compatible with the obvious actions. In particular, the associated characteristic polynomial
$$det(tId -T_0^m)=det(tId -T_{\infty}^m   )$$
coincides to the $m$-th
 Alexander polynomial of $X$ in $U$  for $m<n$.

\end{corollary}

The next result can be regarded as similar to some results in  \cite{CDO}, \cite{Li4} and  \cite{DLi}.

\begin{theorem}  \label{patruth}

Let $g=g_0...g_m=0$ be the equation of the  hypersurface arrangement $W$ in $Z$ and let $F(g)$ be the corresponding
global Milnor fiber given by $g=1$ in the cone $CZ$ over $Z$.
Then 
$$H^j(F(g),\C)=H^j(M_W,\C)$$
 for all $j<n+1$.
In other words, the action of the monodromy on $H^j(F(g),\C)$ is trivial  for all $j<n+1$.

\end{theorem}

Now we turn to the following reformulation of Theorem 4.2 in \cite{M} in our more general setting
described above.

\begin{theorem} \label{doi}

Assume that $d_0$ divides the sum $\sum_{j=1,m}d_j$, say $d d_0=\sum_{j=1,m}d_j$.
Let $ \lambda \in \C^{*}$ be such that $\lambda ^d = 1$ and let $\sigma$ be a non negative integer.
Assume that $ \lambda   $ is not a root of the $q$-th  local Alexander polynomial $\Delta _q (t)_x$ of the hypersurface singularity $(V,x)$ for any $q <n+1 - \sigma$ and any point $x \in W_1$, where $W_1$ is an irreducible component of $W$ different from $W_0$.

Then $ \lambda   $ is not a root of the global Alexander polynomials   $\Delta _q (t)$ associated to $X$ for any
 $q <n+1 - \sigma$.
 \end{theorem}

The  proofs we propose below are based on Theorem 4.2 in \cite{DN} (which relates Alexander modules to the cohomology of a class of rank one local systems
on the complement $M_W$) and on a general idea of getting vanishing results
 via perverse sheaves (based on Artin's vanishing Theorem) introduced in  \cite{CDO}  and developped in  \cite{D}, Chapter 6. We use mainly the notation from  \cite{DN}.

\subsection {Open problem}

It would be interesting to compare the following properties for a polynomial function
$f: \C^{n+1} \to \C$.

\bigskip

\noindent (i) $f$ is $\infty$-transversal;

\bigskip

\noindent (ii) $f$ is h-good as defined in  \cite{DN};

\bigskip

 \noindent (iii) $f$ is $M_0$-tame as defined in Theorem 2.1 in \cite{D3}. 

\bigskip

By taking $X$ to be a hyperplane arrangement with at least two parallel hyperplanes, we see that
(iii) does not imply (i). By taking $X$ to be a central hyperplane arrangement, we see that
(i) does not imply (ii). This example also shows that for an $\infty$-transversal polynomial $f$,
the closure of the nearby fibers $X_t=f^{-1}(t)$ are not necessarily transverse to the hyperplane at infinity, i.e. $f-t$ is not necessarily an $\infty$-transversal polynomial for $t \ne 0$.

On the positive side, note  that  an $\infty$-transversal polynomial $f$  satisfies a {\it partial rational version} of the conditions of being h-good as defined in  \cite{DN}. More precisely, the necessary condition $H_q(T_b,F;\Z)=0$ for all bifurcation points $b$ of $f$ and all $q< n+1$ in the case of an h-good mapping is replaced by the condition  $H_q(T_b,F;\Q)=0$ for all {\it non zero } bifurcation points $b$ of $f$ and all $q< n+1$ (which follows from Theorem 2.10.v in  \cite{DN} and Theorem \ref{treith}, (ii) above.

 A positive answer to the implication (i) implies (iii) would give another proof for Corollary 3.8 in view of Theorem 2.10 v  in  \cite{DN}.

\section{The proofs}

We start with the following easy remark, which proves the first claim of Theorem \ref{treith} in the
case considered by Maxim, i.e. $W_0=H$ is the hyperplane at infinity (for this reason we use here the notation from the beginning of our paper.

\begin{lemma}  \label{zero}
With the above notation, if the closure $V$ of $X$ in $ \PP^{n+1}$ has a positive
dimension singular locus, i.e. $dim V_{sing} >0$, and $H$ is transversal to $V$ except at finitely many points, then 
$$dim V_{sing} =dim (V_{sing} \cap H )+1.$$
 In particular, the singular locus $V_{sing} $ cannot be contained in $H$.
\end{lemma}

Let $f=0$ be a reduced equation for the affine hypersurface $X$. It follows from the above result that $f$ is a primitive polynomial, i.e. its generic fiber $F$ is connected and hence the results in \cite{DN} apply to this situation.

Now we start the proof of  Theorem \ref{treith}. In order to establish the first claim,
note that the closure ${\overline F}$ of $F$ is a general member of the pencil
$$g_1(x)...g_m(x)-tg_0(x)^d=0.$$
As such, it is smooth outside the base locus given by $g_1(x)...g_m(x)=g_0(x)=0$.
A closer look shows that a singular point is located either at a point where at least two of the
polynomials $g_j$ for $1 \leq j  \leq n$ vanish, or at a singular point on one of the hypersurfaces
$W_j$ for $1 \leq j  \leq n$. It follows that $codim Sing({\overline F}) \geq 3$, hence
 ${\overline F}$ is irreducible. This implies that $F$ is connected.

According to  Theorem 4.2 in  \cite{DN}, to prove  the second claim in Theorem \ref{treith}, it is enough to prove the following.

\begin{prop}  \label{unu}

Let $ \lambda \in \C^{*}$ be such that $ \lambda^d \ne 1$, where $d$ is the quotient of $\sum_{j=1,m}d_j$
by $d_0$ .
If $\LL_{\lambda} $ denotes the corresponding local system on $M_W$,
then $H_q(M_W,\LL_{\lambda}  )=0$ for all $q \ne n+1$.

\end{prop}

\proof

First we shall recall the construction of the rank one local system  $\LL_{\lambda} $. Any such local system 
on $M_W$ is given by a homomorphism from $\pi _1(M_W)$ to $\C^*$. To define our local system consider the composition

$$ \pi _1(M_W) \to  \pi _1(M_W') \to  H _1(M_W')= \Z^{m+1}/(d_0,...,d_m) \to \C^*$$
where the first morphism is induced by the inclusion, the second is the passage to the abelianization
and the third one is given by sending the classes $e_0, ...,e_m$ corresponding to the canonical basis of $\Z^{m+1}$ to ${\lambda}^{-d},{\lambda},..., 
{\lambda}$ respectively. For the isomorphism in the middle, see for instance \cite{D0}, p. 102.

\bigskip

It is of course enough to show the vanishing in cohomology, i.e.
 $H^q(M_W,\LL_{\lambda}  )=0$   for all $q \ne n+1$. Let $i:M_W \to  U$ and $j: U \to Z$ be the two inclusions. Then one clearly has
$ \LL_{\lambda}   [n+1] \in Perv(M_W)$ and hence $\F=Ri_*(  \LL_{\lambda}  [n+1]) \in Perv( U   )$,
since the inclusion $i$ is a quasi-finite affine morphism. See for this and the following p. 214 in  \cite{D} for a similar argument.

Our vanishing result will follow from a study of the natural morphism
$$ Rj_!\F \to  Rj_*\F.$$ 
Extend it to a distinguished triangle
$$ Rj_!\F \to  Rj_*\F   \to \G  \to  . $$
Using the long exact sequence of hypercohomology coming from the above triangle, we see exactly as on p.214 in \cite{D} that all we have to show is that
$\HH^k(Z,\G)=0$ for all $k<0$. This vanishing obviously holds if we show that $\G=0$.

This in turn is equivalent to the vanishing of all the local cohomology groups of $ Rj_*\F$, namely $H^m(M_x,\LL_x)=0$ for all $m \in \Z$ and for all points $x \in W_0$.
Here  $M_x = M_W \cap B_x$ for $B_x$ a small open ball at $x$ in $ Z$ and $\LL_x$ is the restriction of the local system $\LL_{\lambda}  $ to $M_x$.

The key observation is that, as already stated above, the action of an oriented elementary loop about the hypersurface $W_0$ in the  local systems  $\LL_{\lambda}  $ and $\LL_x$ corresponds to multiplication by 
$\nu = \lambda ^{-d} \ne 1$.

There are two cases to consider.

\bigskip

\noindent{Case 1.} If $x \in W_0 \setminus V $, then $M_x$ is homotopy equivalent to
 $\C^*$ and the corresponding local system $\LL_{\nu}  $ on  $\C^*$ is defined by multiplication by $\nu$, hence the claimed vanishings are obvious.

\bigskip

\noindent{Case 2.} If  $x \in W_0 \cap V $, then due to the local product structure of stratified sets cut by a transversal, $M_x$ is homotopy equivalent to a product $(B' \setminus (V \cap B')) \times \C^*$, with $B'$ a small open ball centered at $x$ in $W_0$, and the corresponding local system is an external tensor product, the second factor being exactly $\LL_{\nu}    $.  The claimed vanishings follow then from the K\"unneth Theorem, see 4.3.14 \cite{D}.

\bigskip

This ends the proof of  Proposition   \ref{unu}  and of  Theorem \ref{treith}.

\bigskip

Now we turn to the proof of Corollary \ref{top}. The first claim follows from Proposition 6.3.6 and Exercise 4.2.13 in  \cite{D} in conjunction to  Theorem 2.10 v  in  \cite{DN}. In fact, to get the vanishing of
$(R^mf_*\C_U)_0$ one has just to write the exact sequence of the triple $(U,T_0,F)$ and to use the fact that
$H^m(U,\C)=0$ for $m<n+1.$ This vanishing comes from identifying  $U$ to a finite cyclic quotient of the Milnor fiber $F(g_0)$ of the homogeneous isolated complete intersection singularity $CW_0$ defined by
$$ g_0: (CZ,0) \to (\C,0).$$
 For the second claim, one has to use Theorem 2.10.i and Proposition 2.18 in  \cite{DN}.

\bigskip

The proof of  Theorem \ref{patruth} is completely similar to the second part of the proof above and we can safely leave it to the reader after the following remark. Let $D=\sum _{j=0,m}d_j$ and let 
$\alpha $ be a $D$-root of unity, $\alpha \ne 1$. All we have to show is that  $H^q(M_W,\LL_{\alpha}  )=0$   for all $q \ne n+1$, see for instance 6.4.6 in  \cite{D}.

   The action of an oriented elementary loop about the hypersurface $W_0$ in the  local systems  $\LL_{\alpha}  $ and in its restrictions $\LL_x$ as above corresponds to multiplication by 
$\alpha \ne 1$. Therefore the above proof works word for word.

\bigskip

Now we turn to the proof of Theorem  \ref{doi}. We start by the following general remark.

\begin{remark}\label{unurm}

 If $S$ is an $s$-dimensional stratum in a Whitney stratification of $V$ such that $x \in S$ and $W_0$ is transversal to $V$ at $x$, then, due to the local product structure,  the $q$-th reduced local Alexander polynomial $\Delta _q (t)_x$ is the same as that of the  hypersurface singularity $V \cap T$ obtained by cutting
the germ $(V,x)$ by an $(n+1-s)$-dimensional transversal $T$. It follows that these reduced local Alexander polynomials $\Delta _q (t)_x$ are all trivial except for $q \leq n-s$. It is a standard fact that, in the local situation of a hypersurface singularity, the  Alexander polynomials can be defined either from the link or as the characteristic polynomials of the corresponding  the monodromy operators. Indeed, the local Milnor fiber is homotopy equivalent to the corresponding infinite cyclic covering.

\end{remark}

 Let $i:M_W \to  Z \setminus W_1  $ and $j:  Z \setminus W_1      \to Z$ be the two inclusions. Then one  has
$ \LL_{\lambda}    [n+1] \in Perv(M_W)$ and hence $\F=Ri_*( \LL_{\lambda}   [n+1]) \in Perv( Z \setminus W_1    )$,
exactly as above.

Extend now the natural morphism $ Rj_!\F \to  Rj_*\F$ to a distinguished triangle
$$ Rj_!\F \to  Rj_*\F   \to \G .   $$

Applying Theorem 6.4.13 in  \cite{D} to this situation, and recalling the above use of  Theorem 4.2 in \cite{DN}, all we have to check is that $H^m(M_x,\LL_x)=0$ for all points $x \in W_1$  and $m<n+1 - \sigma$.  For $ x \in W_1 \setminus W_0$, this claim is clear by the assumptions made. The case when $x \in W_1 \cap W_0$  can be treated exactly as above, using the product structure, and the fact that
the monodromy of $(W_1,x)$ is essentially the same as that of  $(W_1 \cap W_0 ,x)$, see our remark above.

This completes the proof of  Theorem  \ref{doi}.

\begin{remark}\label{doirm}
Here is an alternative explanation for some of the bounds given in  Theorem 4.2 in \cite{M}.
Assume that $ \lambda $ is a root of the Alexander polynomial 
 $\Delta _i (t)$ for some $i<n+1$. Then it follows from Proposition \ref{doi}
the existence of a point $x \in W_1$ and of an integer $\ell \leq i$ such that
 $\lambda  $ is a root of the local Alexander polynomial $\Delta _{\ell} (t)_x$.
If $x \in S$, with $S$ a stratum of dimension $s$, then by Remark \ref{unurm}, we have $\ell \leq n-s$. This provides half of the bounds in  Theorem 4.2 in \cite{M}. The other half comes from the following remark. Since 
 $\lambda  $ is a root of the Alexander polynomial 
 $\Delta _i (t)$, it follows that $H^i(M_W, \LL_{\lambda}      ) \ne 0$. This implies via an obvious exact sequence that $\HH^{i-n-1}(W_1,\G)\ne 0$. Using the standard spectral sequence to compute this hypercohomology group, we get that some of the groups
$H^p(W_1,\HC^{i-n-1-p}\G)$ are non zero. This can hold only if $p \leq 2 dim (Supp\HC^{i-n-1-p}\G)$. Since $\HC^{i-n-1-p}\G_x = H^{i-p}(M_x,L_x)$ his yields the inequality $p=i- \ell \leq 2s$ in
 Theorem 4.2 in \cite{M}. 
\end{remark}

\begin{remark}\label{treirm}

Let $ \lambda \in \C^{*}$ be such that $ \lambda^d = 1$, where $d$, the quotient of $\sum_{j=1,m}d_j$
by $d_0$, is assumed to be an integer.
Let $\LL_{\lambda} $ denotes the corresponding local system on $M_W$. The fact that the associated monodromy about the divisor $W_0$ is trivial can be restated as follows. Let  $\LL_{\lambda}' $ be the rank one local system on $M_V=Z \setminus V$ associated to ${\lambda}$. Let $j:M_W \to M_V$ be the inclusion. Then
$$\LL_{\lambda}=j^{-1}\LL_{\lambda}'.$$
Let moreover  $\LL_{\lambda}'' $ denote the restriction to  $\LL_{\lambda}' $ to the smooth divisor 
$W_0 \setminus (V \cap W_0)$. Then we have the following Gysin-type long exact sequence
$$ ...  \to  H^q(M_V,\LL_{\lambda}') \to  H^q(M_W,\LL_{\lambda}) \to H^{q-1}(W_0 \setminus (V \cap W_0),\LL_{\lambda}'')  \to  H^{q+1}(M_V,\LL_{\lambda}')    \to    ...$$
exactly as in  \cite{D} , p.222.

 The cohomology groups $ H^*(M_V,\LL_{\lambda}')$ and $ H^*(W_0 \setminus (V \cap W_0),\LL_{\lambda}'')$ being usually simpler to compute than $ H^*(M_W,\LL_{\lambda})$, this exact sequence can give valuable information on the latter cohomology groups.

\end{remark}

\end{document}